\title{Multigraphs (only) satisfy a weak triangle removal lemma}

\author{Asaf Shapira
\thanks{School of Mathematics and College of Computing,
Georgia Institute of Technology, Atlanta GA, 30332. E--mail:
asafico@math.gatech.edu} \and Raphael Yuster\thanks{Department of
Mathematics, University of Haifa, Haifa 31905, Israel. E--mail:
raphy@math.haifa.ac.il} }

\date{}

\documentclass [letterpaper,11pt]{article}
\usepackage{amsfonts}

\newtheorem{theo}{Theorem}

\newcommand{\qed}{\hspace*{\fill} \rule{7pt}{7pt}}

\newcommand{\ignore}[1]{}

\begin{document}
\maketitle

\begin{abstract}

The triangle removal lemma states that a simple graph with
$o(n^3)$ triangles can be made triangle-free by removing $o(n^2)$
edges. It is natural to ask if this widely used result can be
extended to multi-graphs (or equivalently, weighted graphs). 
In this short paper we rule out the
possibility of such an extension by showing that there are
multi-graphs with only $n^{2+o(1)}$ triangles that are still far
from being triangle-free. On the other hand, we show that for some
$g(n)=\omega(1)$, if a multi-graph (or weighted graph) has only $g(n)n^2$
triangles then
it must be close to being triangle-free. The proof relies on
variants of the Ruzsa-Szemer\'edi theorem \cite{RuS}.

\end{abstract}

\section{Introduction}\label{intro}

Motivated by a problem in the theory of extremal hypergraphs,
Ruzsa-Szemer\'edi \cite{RuS} proved the following two theorems.

\begin{theo}[Ruzsa-Szemer\'edi \cite{RuS}]\label{removallemma}
If $G$ is an $n$ vertex graph from which one should remove at least
$\epsilon n^2$ edges in order to destroy all triangles, then $G$
contains at least $f(\epsilon) n^3$ triangles.
\end{theo}

\begin{theo}[Ruzsa-Szemer\'edi \cite{RuS}]\label{RuzSzeBeh}
Suppose $S \subseteq [n]$ is a set of integers containing no 3-term
arithmetic progression. Then there is a graph $G=(V,E)$ with
$|V|=6n$ and $|E|=3n|S|$, whose edges can be (uniquely) partitioned
into $n|S|$ edge disjoint triangles. Furthermore, $G$ contains no
other triangles.
\end{theo}

These two theorems turned out to be two of the most influential
results in extremal combinatorics. First, a simple application of
these two theorems gives a short proof of Roth's Theorem \cite{Roth}
stating that a subset of $[n]$ of size $\epsilon n$ contains a
$3$-term arithmetic progression. The results in \cite{RuS} were
followed by a long line of investigations leading to the recent
hypergraph removal lemmas \cite{Gowers,NRS,RS,Tao}, that also lead
to new proofs of Szemer\'edi's Theorem \cite{Sztheo} and some of its
extensions.

Besides the above applications to additive number theory and
extremal hypergraph theory, which were the original motivation for Theorems
\ref{removallemma} and \ref{RuzSzeBeh}, they also turned out to have
many additional surprising applications. In particular, these
theorems also had applications in extremal combinatorics
\cite{Fure,CEMZ}, in the study of probabilistically checkable proofs
and analysis of linearity tests \cite{HW}, in communication
complexity \cite{PS}, as well as in testing monotonicity
\cite{FLNRRS} and testing graph properties \cite{Alon,AS}.

Theorem \ref{removallemma}, also known as the {\em triangle removal
lemma}, was originally proved for simple graphs, that is, graphs
containing no parallel edges. The proof of Theorem
\ref{removallemma} applies the regularity lemma \cite{Sz}, which can
only handle graphs with constant edge multiplicity\footnote{The edge
multiplicity of a graph is the maximum number of parallel edges
between any pair of vertices.}. In many applications one thus has to
be careful and argue that the graph (or hypergraph) on which one
tries to apply Theorem \ref{removallemma} is indeed simple; see
\cite{S} for one such example. It is thus natural to ask if the
removal lemma also holds for multi-graphs with possibly {\em
unbounded} edge multiplicity. Another way of thinking about this
question is whether the removal lemma holds when the edges of a
graph have {\em arbitrary} weights. Note that if we were to identify
triangles with their vertex sets, then a simple counter example to
such a removal lemma would be to take three vertices, and connect
each pair with $n^2$ edges. In multi-graphs, however, we identify a
triangle with its set of edges\footnote{Note that in simple graphs
there is no difference between identifying a triangle with its edge
set or its vertex set.}, so the above example actually has $n^6$
triangles. We first show that even with this way of counting the
number of triangles, the removal lemma does not hold in
multi-graphs\footnote{It was actually stated, without proof, in some
papers (see, e.g., \cite{Kral2}) that the removal lemma holds only
in simple graphs, although we are not aware of any proof of this
fact. }.

\begin{theo}\label{theoparaupper}
There exists a multi-graph $G$ on $n$ vertices, which contains only
$n^22^{\sqrt{8\log n}}=n^{2+o(1)}$ triangles, and yet one should
remove $n^2$ edges from $G$ in order to make it triangle-free.
\end{theo}

We note that the edge multiplicity of the multi-graph we use in the
proof of Theorem \ref {theoparaupper} is $2^{\sqrt{8\log
n}}=n^{o(1)}$, so we see that the removal lemma fails even when the
edge multiplicity is sub-linear in the size of the graph.

Observe that if we need to remove $n^2$ edges from a graph in order to
make it triangle-free, then it trivially contains at least $ n^2$
triangles. While Theorem \ref{theoparaupper} states that this
trivial lower bound cannot be substantially improved, we can still
ask if a minor improvement is possible. The main motivation is that in some
cases (e.g., the original one in \cite{RuS}) one actually only needs to know
that if a graph is far from being triangle free then it contains asymptotically
more
than $n^2$ triangles. The following theorem answers this question positively.

\begin{theo}\label{theoparalower} If $G$ is an $n$-vertex multi-graph
from which one should remove at least $n^2$ edges in order to
destroy all triangles, then $G$ contains $\omega(n^2)$ triangles.
\end{theo}

We note that because the proof of Theorem \ref{theoparalower}
applies Theorem \ref{removallemma}, the improvement we obtain is
very minor and gives a lower bound of roughly $n^2(\log^*n)^c$ for
some $c>0$ on the number of triangles in the graph. We also remind the reader
that Theorem \ref{theoparalower} can also be stated with respect to weighted
(simple) graphs rather than multigraphs.

\section{The proofs}

For completeness we start with the short proof of Theorem \ref{RuzSzeBeh}.

\paragraph{Proof of Theorem \ref{RuzSzeBeh}:} We define a 3-partite
graph $G$ on vertex sets $A$, $B$ and $C$, of sizes $n$, $2n$ and
$3n$ respectively, where we think of the vertices of the sets $A$,
$B$ and $C$ as representing the sets of integers $[n]$, $[2n]$ and
$[3n]$. For every $1 \leq i \leq n$ and $s \in S$ we put a triangle
$T_{i,s}$ in $G$ containing the vertices $i \in A$, $i+s \in B$ and
$i+2s \in C$. It is easy to see that the above $n|S|$ triangles are
edge disjoint, because every edge determines $i$ and $s$. To see
that $G$ does not contain any more triangles, let us observe that
$G$ can only contain a triangle with one vertex in each set. If the
vertices of this triangle are $a \in A$, $b \in B$ and $c \in C$,
then we must have $b=a+s_1$ for some $s_1 \in S$,
$c=b+s_2=a+s_1+s_2$ for some $s_2 \in S$, and
$a=c-2s_3=a+s_1+s_2-2s_3$ for some $s_3$. This means that
$s_1,s_2,s_3 \in S$ form an arithmetic progression, but because $S$
is free of 3-term arithmetic progressions it must be the case that
$s_1=s_2=s_3$ implying that this triangle is one of the triangles
$T_{i,s}$ defined above. $\qed$

\bigskip

For the proof of Theorems \ref{theoparaupper}
we will need to combine Theorem \ref{RuzSzeBeh}
with the following well known result of Behrend \cite{B} that was recently slightly improved by Elkin \cite{El}.

\begin{theo}[Behrend \cite{B}, Elkin \cite{El}]\label{Behrend}
For every $n$, there exists $S \subseteq [n]$ of size
$n/2^{\sqrt{8\log n}}=n^{1-o(1)}$ containing no 3-term arithmetic
progression.
\end{theo}

\paragraph{Proof of Theorem \ref{theoparaupper}:} Let
$G'$ be the graph of Theorem \ref{RuzSzeBeh} when taking $S
\subseteq [n]$ to be a $3AP$-free set of size $n/2^{\sqrt{8\log n}}$
as guaranteed by Theorem \ref{Behrend}. Let $G$ be the graph
obtained by replacing every edge of $G'$ with $n/|S|=2^{\sqrt{8\log
n}}$ parallel edges. Observe that as $G'$ contains $n|S|$ edge
disjoint triangles, one must remove at least $n|S|$ edges from it in
order to make it triangle-free. As $G$ contains $n/|S|$ parallel
edges for every edge of $G'$ we infer that one must remove $n^2$
edges from $G$ in order to make it triangle-free. Finally, as $G'$
contains only $n|S|$ triangles, we infer that $G$ contains only
$n|S|(n/|S|)^3=n^22^{2\sqrt{8\log n}}$ triangles, as needed. $\qed$

\paragraph{Proof of Theorem \ref{theoparalower}:} Given a
multi-graph $G$, let $T$ be the simple graph on the same vertex set
that contains an edge $(u,v)$ if and only if $G$ has at most
$g^2(n)$ edges connecting $u$ and $v$ for some function
$g(n)=\omega(1)$ to be chosen shortly. Let's first consider the case
that one needs to remove at least $\frac{1}{2g^2(n)}n^2$ edges from
$T$ in order to make it triangle-free. In this case, by Theorem
\ref{removallemma}, we know that $T$ contains at least
$f(\frac{1}{2g^2(n)})n^3$ triangles. Let us now choose a function
$g(n)=\omega(1)$ such that $f(\frac{1}{2g^2(n)})n^3=\omega(n^2)$.
This is clearly possible no matter how fast $f(\epsilon)$ goes to 0
with $\epsilon$. Specifically, given the known bounds on
$f(\epsilon)$ in Theorem \ref{removallemma} (see, e.g., \cite{KS}),
one can take $g(n)=(\log^*n)^c$ for some constant $c>0$. Fixing this
choice of $g(n)$ guarantees that in this case $T$ contains
$\omega(n^2)$ triangles and so $G$ contains at least this many
triangles as well.

So we can assume that we can remove from $T$ a set of edges $E$ of
size $\frac{1}{2g^2(n)}n^2$ and thus make it triangle-free. Let us
now remove from $G$ all the edges connecting pairs of vertices that
are connected by $E$ in $T$. Note that we thus remove from $G$ at
most $g^2(n) \cdot \frac{1}{2g^2(n)}n^2 \leq n^2/2$ edges, hence the
new graph we obtain, let's call it $G'$, has the property that we
should remove at least $n^2/2$ edges from it in order to make it
triangle-free. Furthermore, each edge in $G'$ has multiplicity at
least $g^2(n)$.

Let $T'$ be the simple graph underlining $G'$, that is, the graph on
the same vertex set, with an edge $(u,v)$ if and only if $G'$ has an
edge between $u$ and $v$. Assume first that $T'$ contains at least
$n^2/g(n)$ edges that belong to a triangle. In this case $T'$
contains at least $n^2/3g(n)$ triangles, and as the edge multiplicity
of $G'$ is at least $g^2(n)$ this means that $G'$ contains at least
$n^2g(n)/3$ triangles. As $G'$ is a subgraph of $G$ we infer that
$G$ also contains $n^2g(n)/3$ triangles.

So we can now assume that $T'$ has at most $n^2/g(n)$ edges that
belong to a triangle. Let $E'$ be a set of minimal size whose
removal from $G'$ makes it triangle-free. Let $B$ denote the set of
pairs $(u,v)$ for which $E'$ contains at least one edge connecting
$u$ and $v$, and note that by our assumption on $T'$ we have that
$|B| \leq n^2/g(n)$. For each pair of vertices $(u,v) \in B$ let
$m_{u,v}$ be the number of edges connecting $u$ and $v$ that belong
to $E'$. We claim that for every $(u,v)$ there are at least
$m_{u,v}$ paths of length exactly 2 connecting $u$ and $v$. Indeed,
if $G'$ contains less than $m_{u,v}$ such paths, then we can remove
the $m_{u,v}$ edges connecting $u$ and $v$ from $E'$ and replace
them by one edge from each of the paths of length 2 connecting $u$
and $v$. The new set has fewer edges and it still makes $G'$
triangle-free, which contradicts the minimality of $E'$. We thus
conclude that for every pair $u,v$ the graph $G'$ has at least
$m^2_{u,v}$ triangles containing $u$ and $v$. Recall that $G'$ still
has the property that one should remove at least $n^2/2$ edges from
it in order to make it triangle-free. Therefore we have $\sum
m_{u,v} =|E'| \geq n^2/2$. Combining the above facts, and using
Cauchy-Schwartz, we conclude that the number of triangles in $G'$
(and so also in $G$) is at least
$$
\sum_{(u,v) \in B} m^2_{u,v} \geq \frac{1}{|B|}\left(\sum_{(u,v) \in
B}m_{u,v}\right)^2 \geq \frac{1}{4|B|}n^4 \geq \frac14g(n)n^2\;,
$$
thus completing the proof. $\qed$


\end{document}